\documentclass[graybox]{svmult}

\usepackage{type1cm}        % activate if the above 3 fonts are
\usepackage{makeidx}         % allows index generation
\usepackage{graphicx}        % standard LaTeX graphics tool
\usepackage{multicol}        % used for the two-column index
\usepackage[bottom]{footmisc}% places footnotes at page bottom

\usepackage{newtxtext}       % 
\usepackage[varvw]{newtxmath}       % selects Times Roman as basic font

\usepackage{xcolor}
\usepackage{url}
\usepackage{tabu}
\usepackage{bm}
\usepackage{algorithm}
\usepackage{algorithmic}

\usepackage{amsmath,amsfonts}

\usepackage{todonotes}

\usepackage{hyperref}
\usepackage{cleveref}

\definecolor{dred}{rgb}{0.8,0,0}
\def\CRR{\color{black}}

\makeindex             % used for the subject index
\begin{document}

\title*{Finite basis physics-informed neural networks as a Schwarz domain decomposition method}
% Domain decomposition training strategies for physics-informed neural networks

\titlerunning{FBPINNs as a Schwarz method}

\author{Victorita Dolean, Alexander Heinlein, Siddhartha Mishra, and Ben Moseley}
\institute{Victorita Dolean \at University of Strathclyde, 26, Richmond Street, G1 1XH Glasgow, UK and University Côte d'Azur, CNRS, LJAD, France. \email{work@victoritadolean.com}
\and Alexander Heinlein \at Delft University of Technology, Faculty of Electrical  Engineering  Mathematics \& Computer Science, Delft Institute of Applied Mathematics, Mekelweg 4, 2628 CD Delft, Netherlands.  \email{a.heinlein@tudelft.nl}
\and 
Siddhartha Mishra and Ben Moseley \at ETH Z\"urich, Computational and Applied Mathematics Laboratory / ETH AI Center, R\"amistrasse 101, 8092 Z\"urich, Switzerland.}

\maketitle

\abstract*{Physics-informed neural networks (PINNs) \cite{Lagaris1998, raissi2019physics} are an approach for solving boundary value problems based on differential equations (PDEs). The key idea of PINNs is to use a neural network to approximate the solution to the PDE and to incorporate the residual of the PDE as well as boundary conditions into its loss function when training it. This provides a simple and mesh-free approach for solving problems relating to PDEs. However, a key limitation of PINNs is their lack of accuracy and efficiency when solving problems with larger domains and more complex, multi-scale solutions.
In a more recent approach, finite basis physics-informed neural networks (FBPINNs) \cite{moseley2021finite} use ideas from domain decomposition to accelerate the learning process of PINNs and improve their accuracy. In this work, we show how Schwarz-like additive, multiplicative, and hybrid iteration methods for training FBPINNs can be developed. We present numerical experiments on the influence of these different training strategies on convergence and accuracy. Furthermore, we propose and evaluate a preliminary implementation of coarse space correction for FBPINNs.}

\vspace{-2cm}

\section{Introduction}

\vspace{-0.2cm}

The success and advancement of machine learning (ML) in fields such as image recognition and natural language processing has lead to the development of novel methods for the solution of problems in physics and engineering. However, algorithms developed in traditional fields of ML usually require a large amount of data, which are difficult to obtain from measurements and/or traditional numerical simulations. Furthermore, such algorithms can be difficult to interpret and can struggle to {\CRR generalize}. To overcome these issues, a new research paradigm has emerged, known as scientific machine learning (SciML) \cite{Baker2019, Moseley2022}{\CRR,} which aims to more tightly combine ML with scientific principles to provide more powerful algorithms. 

%By combining physics and data, SciML has provided promising new approaches for computing observables of complex numerical simulations, solving high-dimensional parametrised PDEs, carrying out data assimilation, solving inverse problems, and discovering incomplete physics. 

One such approach are physics-informed neural networks (PINNs) \cite{Lagaris1998, raissi2019physics}, which are designed to approximate the solution to the boundary value problem
\begin{equation}
\begin{array}{rcll}
\mathcal{N}[u]({\mathbf{x}}) &=& f(\mathbf{x}), \text{ } & \mathbf{x} \in \Omega \subset \mathbb{R}^{d}, \\
\mathcal{B}_k[u] ({\mathbf{x}}) &=& g_k(\mathbf{x}), \text{ } & \mathbf{x}  \in \Gamma_{k} \subset \partial\Omega 
\end{array}
\label{eq:label1}
\end{equation}
where $\mathcal{N}[u](\mathbf{x}) $ is a differential operator, $u$ is the solution and $\mathcal{B}_k(\cdot)$ is a set of boundary conditions, such that the solution $u$ is uniquely determined. Note that boundary conditions are to be understood in a broad sense and the $\mathbf{x}$ variable can also include time. In particular, we do not distinguish between initial and boundary conditions. 

The approximation to the solution of \eqref{eq:label1} is given by a neural network $u(\mathbf{x},\bm{\theta})$ (for the sake of simplicity we use the same notation for the solution of the PDE and the neural network) where $\bm{\theta}$ is a vector of all the parameters of the neural network (i.e., its weights and biases). The network is trained via the loss function 
\begin{align}
   \mathcal{L}(\bm{\theta})= \underbrace{\frac{\lambda_{I}}{N_I} \sum_{i=1}^{N_I} (\mathcal{N}[u](\mathbf{x}_{i},\bm{\theta}) - f(\mathbf{x}_{i}))^{2}}_{\mathcal{L}_{\text{PDE}}}
     + \underbrace{\sum_{k=1}^{N_k} \frac{\lambda_{B}^{k}}{N_{B}^{k}} \sum_{j=1}^{N_{B}^{k}}(\mathcal{B}_k[u](\mathbf{x}_{j}^{k},\bm{\theta})  - g_k(\mathbf{x}_{j}^{k}))^2}_{\mathcal{L}_{\text{BC}}}. \label{eq:loss_pinn}
\end{align}
Here, $\{\mathbf{x}_i\}_{i=1}^{N_I}$ is a set of collocation points sampled in the interior of the domain, $\{\mathbf{x}_{j}^{k}\}_{j=1}^{N_{B}^{k}}$ is a set of points sampled along each boundary condition, and $\lambda_I$ and $\lambda_{B}^{k}$ are well-chosen scalar hyperparameters which ensure that the terms in the loss function are well balanced. Intuitively, one can see that the PDE loss tries to ensure that {\CRR the} solution learned by the network obeys the underlying PDE whilst the boundary loss tries to ensure it obeys the boundary conditions.

In practice, the presence of the boundary loss in \cref{eq:loss_pinn} often slows down training as it can compete with the PDE term \cite{wang2022and}. In a slightly different formulation, boundary conditions can instead be enforced exactly as hard constraints by using the neural network as part of a solution ansatz $\mathcal{C} u$ where ${\cal C}$ is a constraining operator which enforces that the solution explicitly satisfies the boundary conditions \cite{Lagaris1998, moseley2021finite}. This turns the {\CRR optimization} problem into an unconstrained one, and only the PDE loss from \cref{eq:loss_pinn} is required to train the PINN. For example, suppose we want to enforce that $u(0)=0$ when solving a one-dimensional ODE, then the ansatz and constraining operator can be chosen as
$[{\cal C} u] (x,\bm{\theta}) = \tanh(x) u(x,\bm{\theta}).$
The rationale behind this is that the function $\tanh(x) $ is null at $0$, forcing the boundary condition to be obeyed, but non-zero away from $0$, allowing the network to learn the solution away from the boundary condition. 

%Note there is no unique way to reinforce in a hard way the boundary condition and this may not be possible for more complex boundary conditions, i.e. it is problem dependent.

%\begin{equation}
%\label{eq:loss_hard}
%   \mathcal{L}(\bm{\theta}) =\frac{1}{N_I} \sum_{i=1}^{N_I} (\mathcal{N}%[\mathcal{C} u](\mathbf{x}_{i},\bm{\theta}) - f(\mathbf{x}_{i}))^{2}.
%\end{equation}

Whilst PINNs have proven to be successful for solving many different types of differential equations, they often struggle to scale to problems with larger domains and more complex, multi-scale solutions \cite{moseley2021finite, Wang2021d}. This is in part due to the spectral bias of neural networks \cite{pmlr-v97-rahaman19a} (their tendency to learn higher frequencies much slower than lower frequencies), and the increasing size of the underlying PINN {\CRR optimization} function. One way to alleviate these scaling issues is to combine PINNs with a domain decomposition method (DDM); by taking a divide-and-conquer approach, one hopes that the large, global {\CRR optimization} problem can be turned into a series of smaller and easier {\CRR localized} problems. In particular, \cite{moseley2021finite} proposed finite basis physics-informed neural networks (FBPINNs) where the global PINN is partitioned into many local networks that are trained to approximate the solution on an overlapping domain decomposition. Related approaches are the deep domain decomposition method (D3M)~\cite{li2019d3m} and the deep-learning-based domain decomposition method (DeepDDM)~\cite{li2020deep}, which combine overlapping Schwarz domain decomposition methods with a PINN-based discretization. Other earlier works on the use of machine learning and domain decomposition methods include the prediction of the geometrical location of constraints in adaptive FETI-DP and BDDC methods; see~\cite{Heinlein:2019:MLA}. For an overview of the combination of domain decomposition methods and machine learning, see~\cite{Heinlein:2021:CML2}. 

In this work, we build upon FBPINNs by showing how Schwarz-like additive, multiplicative and hybrid iterative training strategies for FBPINNs can be developed. We present numerical experiments on the influence of these training strategies on convergence and accuracy. We propose and evaluate a preliminary implementation of a coarse space correction for FBPINNs, to further improve their efficiency. 

\vspace{-0.5cm}

\section{Finite basis physics-informed neural networks (FBPINNs)}

\vspace{-0.2cm}

First we briefly present the FBPINN method introduced by~\cite{moseley2021finite} from a DDM perspective. The FBPINN method can be seen as a network architecture that allows for a localization of the network training. Therefore, let us consider a set of collocation points $X = \{\mathbf{x}_i\}_{i=1}^N$  in the global domain $\Omega$ and a decomposition into overlapping domains $\Omega = \cup_{j=1}^J \Omega_j$ inducing a decomposition into subsets of collocation points $X_j = \{\mathbf{x}_i^j\}_{i=1}^{N_j}$, $j=1,..., J$. As usual in overlapping Schwarz methods, $X= \cup_{j=1}^J X_j$ is not disjoint. For each subdomain $\Omega_j$, we denote ${\cal N}_j$ the index set of {\CRR neighboring} subdomains, $\Omega_j^{\circ} = \cup_{l=1}^{{\cal N}_j} \Omega_l \cap \Omega_j$ the overlapping subset of $\Omega_j$, and $\Omega_j^{int}  = \Omega_j \setminus \Omega_j^{\circ}$, the interior part of the domain; let $X_j^\circ$ and $X_j^{int}$ be the {\CRR corresponding} sets of collocation points, and $X^\circ = \cup_{j=1}^J X_j^\circ$ and $X^{int} = \cup_{j=1}^J X_j^{int}$.
\begin{figure}[t]
  \centering
  \includegraphics[width=0.49\textwidth]{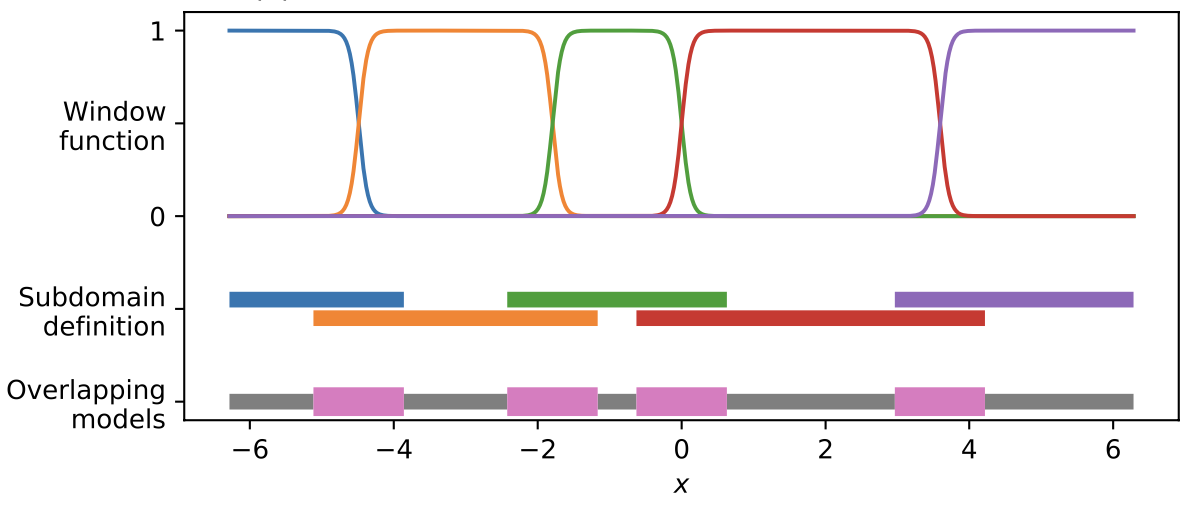}
  \includegraphics[width=0.49\textwidth]{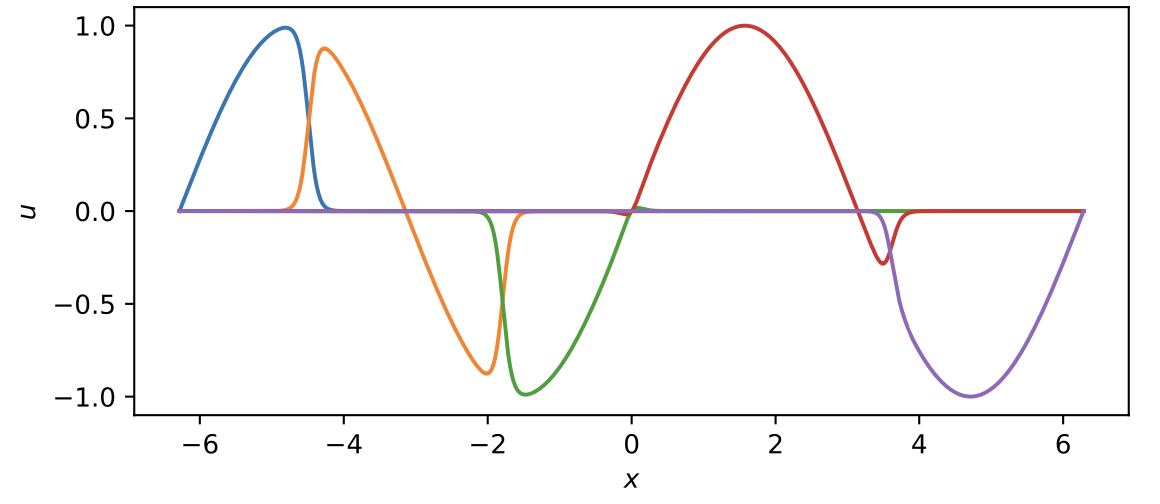}
  \caption{Local FBPINN subdomains and window functions $w_j$ (left), local solutions $u_j$ (right)}
  \label{fig:localnets}
\end{figure}
We now define the global network $u$ as the sum of local networks $u_j(\mathbf{x},\bm{\theta}_j)$ weighted by window functions $\omega_j$:
$
u = \sum_{j, \mathbf{x}_{i} \in \Omega_j} \omega_j u_j.
$
Here, the local networks have individual network parameters $\bm{\theta}_j$, and of course, they could simply be evaluated everywhere in $\mathbb{R}^d$. In order to restrict them to their corresponding overlapping subdomains, we multiply them with the window functions, which have the {\CRR properties} ${\rm supp}(\omega_i) \subset \Omega_i$ {\CRR and $\Omega \subset \cup_{j=1}^J {\rm supp}(\omega_i)$; the specific definition of $\omega_i$ employed here can be found in~\cite[eq.~(14)]{moseley2021finite}.} See~\cref{fig:localnets} for a graphical representation of the overlapping subdomains, their overlapping and interior sets, window functions, and local solutions for a simple one-dimensional example. If we insert {\CRR the expression for $u$ into \cref{eq:loss_pinn}}, we see that {\CRR the loss function} can be written as:
\begin{equation}
\label{eq:loss}
\begin{aligned}
   \mathcal{L} (\bm{\theta}_1,\ldots,\bm{\theta}_J) 
   = & 
	\frac{1}{N} \sum_{i=1}^N \left( \mathcal{N} [\mathcal{C} \sum_{j, \mathbf{x}_{i} \in X_j} \omega_j u_j](\mathbf{x}_{i},\bm{\theta}_j) - f(\mathbf{x}_{i})\right)^{2}.
\end{aligned}
\end{equation}
 The contribution to the global loss function can be split into a part coming from interior points and another one from points in the overlap, respectively:
\begin{equation}
    \begin{aligned}
\mathcal{L} (\bm{\theta}_1,\ldots,\bm{\theta}_J) 	= & 
	\underbrace{\frac{1}{N} \sum_{\mathbf{x} \in X^{int}} \left( \mathcal{N} [\mathcal{C} \sum_{l, \mathbf{x} \in X_l} \omega_l u_l](\mathbf{x},\bm{\theta}_l) - f(\mathbf{x})\right)^{2}}_{=: \mathcal{L}^{int} (\bm{\theta}_1,\ldots,\bm{\theta}_J)} \\
	& 
	+ 
 \frac{1}{N} 
%		\sum_{\Gamma_j^k \subset \Omega_j} \sum_{\mathbf{x}_{i} \in \Gamma_j^k} 
		\sum_{\mathbf{x} \in X^\circ}
		\left( \mathcal{N} [\mathcal{C} \sum_{l, \mathbf{x} \in X_l} \omega_l u_l](\mathbf{x},\bm{\theta}_l) - f(\mathbf{x})\right)^{2}.
\end{aligned}
\end{equation}
Note also that, since $X_i^{int} \cap X_j^{int} = \emptyset$ for $i \neq j$, the interior contribution can be simplified as follows:
$$
	\mathcal{L}^{int} (\bm{\theta}_1,\ldots,\bm{\theta}_J)
	=
	\frac{1}{N} \sum_{j = 1}^J \sum_{\mathbf{x}_{i} \in X_j^{int}} \left( \mathcal{N} [\mathcal{C} \omega_j u_j](\mathbf{x}_{i},\bm{\theta}_j) - f(\mathbf{x}_{i})\right)^{2}.
$$

In \cite{moseley2021finite}, the authors introduce the notion of {\it scheduling} which is related to the degree of parallelism one can consider in Schwarz domain decomposition methods. For example in the well-know alternating Schwarz method, local solves take place sequentially, in an alternating manner, with data being exchanged at the interfaces. In the case of the parallel Schwarz method, local solutions are computed simultaneously, but subdomains only have access to interface data at the previous iteration. As is well-known in DDMs, the alternating method convergences in fewer iterations than the parallel method, whereas the second methods allows the concurrent computation of the local solutions; hence, the parallel Schwarz method is often more efficient in a parallel implementation. 

In the case of many subdomains, one can define a so-called coloring strategy, i.e., subdomains with the same color are computed in parallel and different colors are processed sequentially. Here, we will consider any possible coloring scheme, allowing for arbitrary combinations of additive and multiplicative coupling. In particular, let us split the set of subdomain indices as follows $\{1,...,J\} = {\cal A} \cup {\cal I}$, such that subdomains $\Omega_j,\, j\in {\cal A}$ are allocated the same `color' which is different than those of the subdomains $\Omega_j,\, j\in {\cal I}$. In the case of training FBPINNs, the notion of coloring is replaced by that of  {\it scheduling}, that is, subdomains indexed in ${\cal A}$ are considered to be {\it active} at a given iteration and those indexed in ${\cal I}$ are {\it inactive}. The case when ${\cal I} =\emptyset$ corresponds to the fully parallel Schwarz method, whereas the case where only one subdomain is {\it active} at a time corresponds to a fully alternating Schwarz iteration. Denoting a subdomain $\Omega_j$ as {\it inactive} corresponds to fixing ${\bm \theta}_j$ during the optimization of $\mathcal{L} (\bm{\theta}_1,\ldots,\bm{\theta}_J)$.

The FBPINN training algorithm follows the `coloring' strategy described above. Let us denote by $\bm{\theta}_j^k$ the parameter values at the $k$-th training step, and to simplify the presentation, we focus on the case of a first order gradient-based optimizer. If we start from an initial guess $\bm{\theta}_j^0$, then the training step for each subdomain is given by \cref{alg:training}.
\begin{algorithm}
\caption{FBPINN training step for each subdomain}
\begin{algorithmic}
\IF {$j \in {\cal A}$ ($\Omega_j$ is an active domain)}
\STATE \textcolor{magenta}{Perform $p$ iterations of gradient descent on $\bm{\theta}_j^k$ ($\bm{\theta}_i^{k}$ where $i\ne j$ are kept fixed): 
\vspace{-0.2cm}
$$
\bm{\theta}_j^{k+l} = \bm{\theta}_j^{k+l-1} - \lambda \nabla_{\bm{\theta}_j} \mathcal{L}(\bm{\theta}_1^{k},...,\bm{\theta}_{j-1}^{k},\bm{\theta}_j^{k+l-1},\bm{\theta}_{j+1}^{k},...,\bm{\theta}_J^{k}), l=1,..,p.
$$
}
\vspace{-0.3cm}
\STATE \textcolor{orange}{Update the solution in the overlapping regions (communicate with neighbours):
$$
\forall \mathbf{x}\in \Omega_j^{\circ},\, u(\mathbf{x},\bm{\theta}_j^{k+p}) \leftarrow  \sum_{l, \mathbf{x} \in \Omega_l} \omega_l u_l (\mathbf{x},\bm{\theta}_l^{k+p}).
$$
\vspace{-0.2cm}
}

\ENDIF 
\end{algorithmic}
\label{alg:training}
\end{algorithm}
Once all active subdomains have completed one training step, the set ${\cal A}$ and ${\cal I}$ are updated. This whole procedure is repeated until any stopping criterion, such as a maximum number of iterations or a tolerance for the loss, is met.

\begin{figure}[t]
  \centering
  \includegraphics[width=0.42\textwidth]{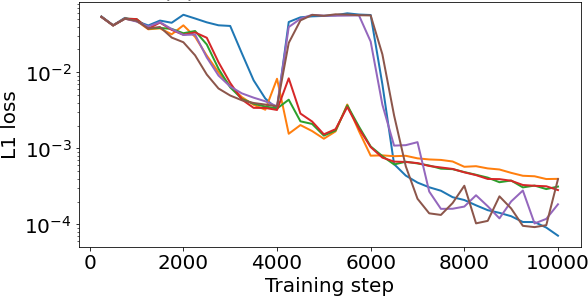}
  \includegraphics[width=0.42\textwidth]{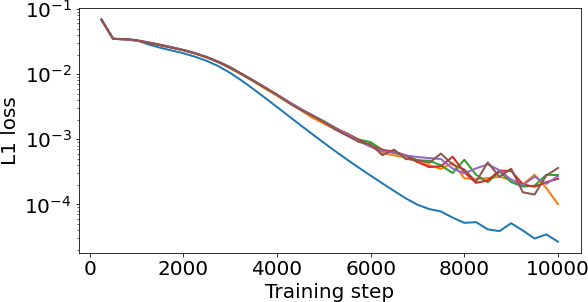} \\[2mm]
   \includegraphics[width=0.42\textwidth]{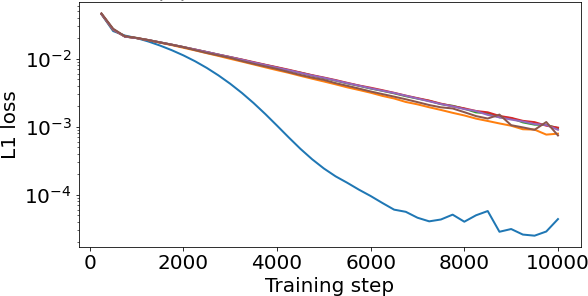}
  \begin{minipage}[c]{0.42\textwidth}
    
    \vspace{-30mm}
  
  \begin{center}
      \includegraphics[width=0.8\textwidth]{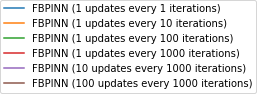}
  \end{center}
  
  \end{minipage}
  \caption{FBPINN convergence for decompositions into $8$ (upper left), $16$ (upper right) and $32$ (lower left) subdomains and different $p$ values. Each network has 2 layers and 16 hidden units per layer. A total of 3,000 collocation points regularly sampled over the domain are used. For each decomposition, subdomains are regularly spaced, with an overlap of 70\% of the subdomain width.
  % batch\_size = 3000 and I have subdomain\_xs = get\_subdomain\_xs([np.array([2,]*16)], [2*np.pi]) and width = 0.7}
  }
  \label{fig:flextrain}
\end{figure}

Let us note that:
\begin{itemize}
\item The gradient updates can be performed in \textcolor{magenta}{parallel} and are fully localized even if the loss function is global; only in the update step are network solutions and network gradients \textcolor{orange}{transferred} between neighboring subdomains.
\item It is not necessary to perform communication in the overlaps (here in \textcolor{orange}{orange}) at every iteration of gradient descent, but rather every $p$ iterations for a better computational efficiency. The overall convergence can also be affected; cf.~\cref{fig:flextrain}. 
\item 
%Unlike classical domain decomposition methods, we always compute gradient updates with respect to the global loss function. Therefore, no local problems are solved in each iteration step.
{\CRR Unlike in classical domain decomposition methods, in our approach, the global problem is not decomposed into local problems, which can be solved independently. Instead, we always compute gradient updates with respect to the global loss function, and the domain decomposition and hence the localization enters through the window functions in the definition of the architecture of global network.}
\end{itemize}

To illustrate the behavior of this algorithm we will consider a scaling study for its flexible training strategy. In particular, we fix the number of global collocation points and investigate the influence of changing the number of subdomains and the value of $p$ on convergence when solving the simple 1D ODE $\frac{du}{dx}=\cos \omega x$, $u(0)=0${\CRR, with $\omega = 15$.}
All subdomains are kept active all of the time, and all other FBPINN design choices are kept the same, including window function and local network architecture per subdomain. {\CRR We only consider the case of relatively large overlap of 70\% of the subdomain width, but the results are qualitatively the same for other sizes of overlap. As discussed in~\cite{moseley2021finite} and as in classical overlapping Schwarz methods, performance generally improves when increasing the size of the overlap; a systematic investigation is still open.}

In Figure \ref{fig:flextrain}, we display the convergence of the loss function when the communication between subdomains takes place every $p\in \{1,10,100,1000\}$ epochs. We observe that the case of $8$ subdomains is rather special since convergence appears rather unstable and there is no option that performs clearly best. As we increase the number of subdomains to $16$ and $32$ we observe an expected behavior, that is, the convergence rate improves if we communicate solutions and gradients in the overlaps every iteration. Moreover, when increasing the number of subdomains, naturally the global training performs less well, which is well known in domain decomposition as lack of scalability; we observe this behavior for all values of $p$. This is expected because the method above corresponds to a one-level method (meaning only neighboring subdomains communicate and there is no global exchange of information). Surprisingly, we do not see any clear difference in the convergence depending on $p$, that is, depending on how often we communicate. Since, in a parallel setting, it is computationally more efficient to communicate less, the results seems to indicate that, if we do not communicate in each step, it is beneficial to communicate as little as possible.

% there is no fundamental difference in the result of the global training when one communicates every $\{10,100,1000\}$ epochs. The second important remark is that when increasing the number of subdomains from $16$ to $32$ for example, naturally the global training performs less well, which is well known in domain decomposition as lack of scalability. This is also expected, as the method we defined corresponds to one-level method (meaning only {\CRR neighboring} subdomains communicate and there is no exchange of global information).

\vspace{-0.5cm}

\section{Coarse correction}

\vspace{-0.2cm}

Coarse spaces are instrumental in DDMs, as they ensure the robustness of a given method with respect to the number of subdomains as well as other problem-specific parameters, such as physical properties like frequency for wave problems or conductivity for diffusion type problems. Coarse spaces are often defined based on geometrical information (like a coarser mesh) but more sophisticated coarse spaces can be constructed using spectral information of underlying local problems. When training PINNs, it {\CRR is} not immediately clear how to define a coarse space, that is, a coarse network model, nor how to choose the number of collocation points and parameters of the coarse model. In what follows, we propose and evaluate a preliminary implementation of a coarse space correction for FBPINNs. 

In particular, we exploit the spectral bias of neural networks in order to build a coarse correction. This is the well-studied phenomenon that they tend to learn higher frequencies much slower than lower frequencies~\cite{pmlr-v97-rahaman19a}, and similar effects are observed for PINNs \cite{moseley2021finite, Wang2021d}. Indeed, this effect is what motivated the use of domain decomposition in FBPINNs. More precisely, we first train a small but global network for enough epochs to learn the low frequency component of the solution{\CRR; in particular, we employ a coarse network with the same architecture as a single local network.} Then, local subdomains are added to approximate missing higher frequency components. The resulting FBPINN solution is given by $
u = u_g + \sum_{j, \mathbf{x}_{i} \in \Omega_j} \omega_j u_j$, where $u_g(\mathbf{x},\bm{\theta}_g)$ is the coarse network and $u_j(\mathbf{x},\bm{\theta}_j)$ are the local networks. Because of spectral bias, low frequencies are first {\CRR learned} by the coarse network, and a relatively small network is sufficient to approximate the low frequencies. Then the local networks only need to learn the remaining higher frequencies. Since the local models only have to learn a local part of the solution, relatively small local network models are also sufficient.

% \begin{algorithm}
% \caption{Coarse correction to FBPINN training}
% \begin{algorithmic}

% \STATE Perform $M$ preliminary iterations (epochs) to train a global, coarse PINN model.
% % on the global problem with the global loss function.

% %\STATE \textcolor{orange}{Update the gradients in the overlap}.
% \STATE Use the current iterate (coarse) $\bm{\theta}_1^M,\bm{\theta}_2^M,...,\bm{\theta}_J^M$ as an initial guess in the FBPINN parallel algorithm. (Algorithm \ref{alg:training})
 
% \end{algorithmic}
% \label{alg:coarse}
% \end{algorithm}

\begin{figure}[t]
  \centering
  \includegraphics[width=0.87\textwidth]{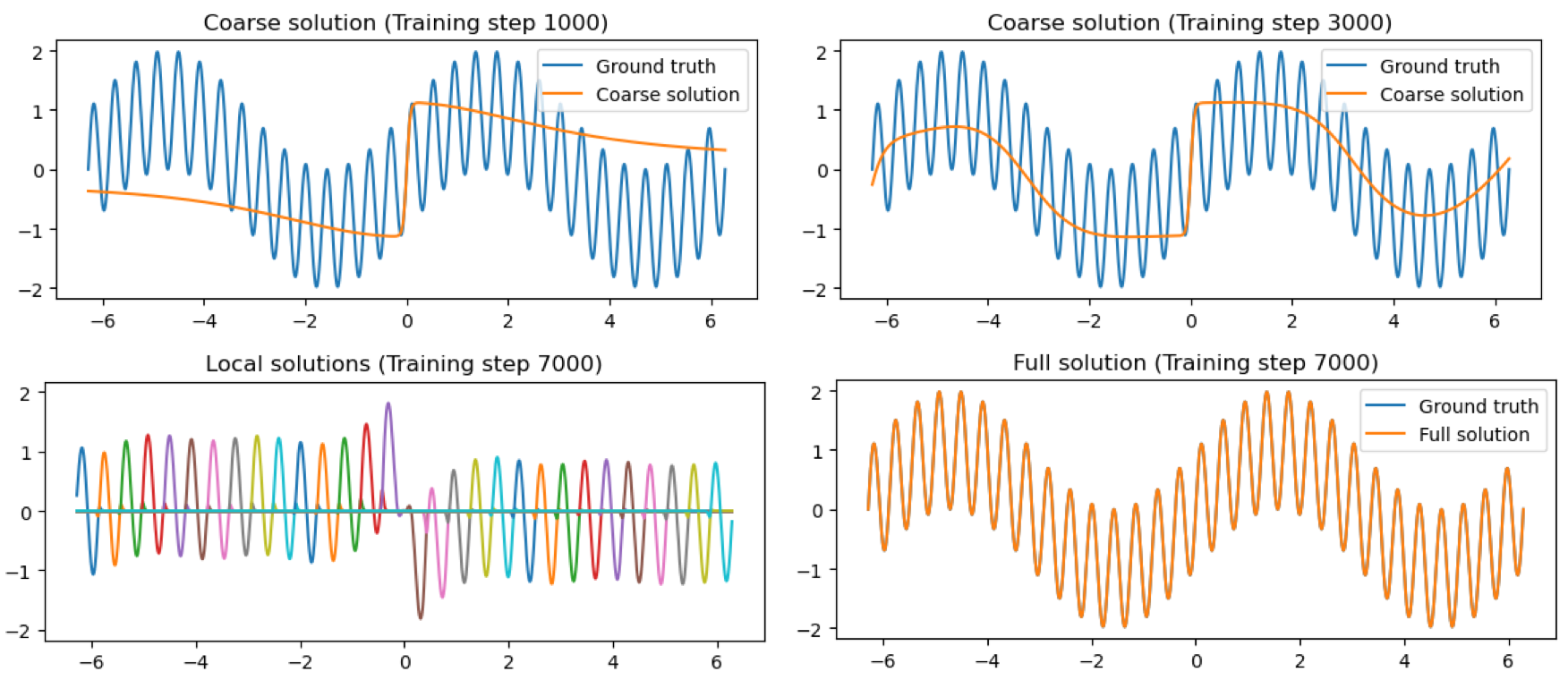}
  \caption{Coarse correction for FBPINNs. Each subdomain network and the coarse network has 2 layers and 16 hidden units per layer. A total of 500 collocation points regularly sampled over the domain are used to train the coarse network, and 3,000 for the local networks. For the local networks, the subdomains are regularly spaced, with an overlap of 70\% of the subdomain width.}
  \label{fig:coarse}
\end{figure}

We will apply these ideas on the simple 1D ODE, $\frac{du}{dx}=\omega_1 \cos(\omega_1 x)+ \omega_2 \cos(\omega_2 x)$, $u(0)=0$, where two frequencies are present in the solution, $u(x) = \sin(\omega_1 x)+ \sin(\omega_2 x)$. For our test case, we choose $\omega_1=1$ as a lower frequency and $\omega_2=15$ as the higher frequency and we decompose the global domain into $30$ overlapping subdomains; see~\cref{fig:coarse}. {\CRR We note, as shown in~\cite{moseley2021finite} for an ODE with a single high frequency, solving such a problem with a single PINN requires a high network complexity and large number of iterations.} First, we train the global {\CRR coarse} network, $u_g$, until the lower frequency is {\CRR learned}. We illustrate this progressive process in \cref{fig:coarse} where we see that we need roughly 3,000 epochs to identify the lower frequency. Here, we have chosen the number of epochs by hand {\CRR based on the accuracy of the coarse solution}, but in the future, we will work on automating the training of the coarse network. Then, the coarse network is fixed and the local networks are trained to approximate the remaining component of the solution, with all local networks kept active at each training step. As can be seen in~\cref{fig:coarse}, using our proposed approach, the coarse network approximates the coarse component of the solution, and the local subdomain networks approximate the high frequency components on the local subdomains.

\vspace{-0.6cm}

\section{Conclusions}

\vspace{-0.2cm}

In this work, we provide first insights on how to incorporate techniques from classical Schwarz domain decomposition methods into the FBPINN method. We show that its algorithmic components can be translated in the language of domain decomposition methods, and the well-established notions of additive, multiplicative and, hybrid Schwarz iterations can be identified through the notion of the flexible scheduling strategies introduced in~\cite{moseley2021finite}. Finally, we start exploring the notion of coarse space for FBPINNs. In particular, we train a coarse network to approximate the low frequency components of the solution and then continue by training local networks to approximate the remaining high frequency components. These ideas can be extended in a straightforward way to other, more complex boundary value problems.

\vspace{-0.6cm}

\bibliographystyle{spmpsci}
\bibliography{ref}

\begin{thebibliography}{10}
\providecommand{\url}[1]{{#1}}
\providecommand{\urlprefix}{URL }
\expandafter\ifx\csname urlstyle\endcsname\relax
  \providecommand{\doi}[1]{DOI~\discretionary{}{}{}#1}\else
  \providecommand{\doi}{DOI~\discretionary{}{}{}\begingroup
  \urlstyle{rm}\Url}\fi

\bibitem{Baker2019}
Baker, N., Alexander, F., Bremer, T., Hagberg, A., Kevrekidis, Y., Najm, H.,
  Parashar, M., Patra, A., Sethian, J., Wild, S., Willcox, K., Lee, S.:
  {Workshop Report on Basic Research Needs for Scientific Machine Learning:
  Core Technologies for Artificial Intelligence}.
\newblock Tech. rep., USDOE Office of Science (SC) (United States) (2019)

\bibitem{Heinlein:2019:MLA}
Heinlein, A., Klawonn, A., Lanser, M., Weber, J.: Machine learning in adaptive
  domain decomposition methods - predicting the geometric location of
  constraints.
\newblock SIAM Journal on Scientific Computing \textbf{41}(6), A3887--A3912
  (2019)

\bibitem{Heinlein:2021:CML2}
Heinlein, A., Klawonn, A., Lanser, M., Weber, J.: Combining machine learning
  and domain decomposition methods for the solution of partial diﬀerential
  equations – a review.
\newblock GAMM-Mitteilungen \textbf{44}(1), e202100001 (2021)

\bibitem{Lagaris1998}
Lagaris, I.E., Likas, A., Fotiadis, D.I.: {Artificial neural networks for
  solving ordinary and partial differential equations}.
\newblock IEEE Transactions on Neural Networks \textbf{9}(5), 987--1000 (1998)

\bibitem{li2019d3m}
Li, K., Tang, K., Wu, T., Liao, Q.: D3m: A deep domain decomposition method for
  partial differential equations.
\newblock IEEE Access \textbf{8}, 5283--5294 (2019)

\bibitem{li2020deep}
Li, W., Xiang, X., Xu, Y.: Deep domain decomposition method: Elliptic problems.
\newblock In: Mathematical and Scientific Machine Learning, pp. 269--286. PMLR
  (2020)

\bibitem{Moseley2022}
Moseley, B.: {Physics-informed machine learning: from concepts to real-world
  applications}.
\newblock Ph.D. thesis, University of Oxford (2022)

\bibitem{moseley2021finite}
Moseley, B., Markham, A., Nissen-Meyer, T.: Finite basis physics-informed
  neural networks (fbpinns): a scalable domain decomposition approach for
  solving differential equations.
\newblock arXiv preprint arXiv:2107.07871  (2021)

\bibitem{pmlr-v97-rahaman19a}
Rahaman, N., Baratin, A., Arpit, D., Draxler, F., Lin, M., Hamprecht, F.,
  Bengio, Y., Courville, A.: On the spectral bias of neural networks.
\newblock In: K.~Chaudhuri, R.~Salakhutdinov (eds.) Proceedings of the 36th
  International Conference on Machine Learning, \emph{Proceedings of Machine
  Learning Research}, vol.~97, pp. 5301--5310. PMLR (2019)

\bibitem{raissi2019physics}
Raissi, M., Perdikaris, P., Karniadakis, G.E.: Physics-informed neural
  networks: A deep learning framework for solving forward and inverse problems
  involving nonlinear partial differential equations.
\newblock Journal of Computational physics \textbf{378}, 686--707 (2019)

\bibitem{Wang2021d}
Wang, S., Wang, H., Perdikaris, P.: {On the eigenvector bias of Fourier feature
  networks: From regression to solving multi-scale PDEs with physics-informed
  neural networks}.
\newblock Computer Methods in Applied Mechanics and Engineering \textbf{384},
  113938 (2021)

\bibitem{wang2022and}
Wang, S., Yu, X., Perdikaris, P.: When and why pinns fail to train: A neural
  tangent kernel perspective.
\newblock Journal of Computational Physics \textbf{449}, 110768 (2022)

\end{thebibliography}

\end{document}